\documentclass[onecolumn,11pt]{article}
\usepackage[top=1in, bottom=1in, left=1in, right=1in]{geometry}
\usepackage{amsmath,amsfonts,amscd,amssymb}
\usepackage{graphicx}
\usepackage{overpic}
\usepackage{cancel}
\usepackage{rotating}
\usepackage{url}
\usepackage{caption}
\usepackage{color}
\usepackage{rotating}
\usepackage{multirow}
\usepackage{wrapfig}
\usepackage{mathtools}
\usepackage{subeqnarray}
\usepackage{setspace}
\usepackage{palatino} 
\graphicspath{{figures/}}

\usepackage[bottom,flushmargin,hang,multiple]{footmisc}
\usepackage{lipsum}
\newcommand\blfootnote[1]{%
  \begingroup
  \renewcommand\thefootnote{}\footnote{#1}%
  \addtocounter{footnote}{-1}%
  \endgroup
}

\newcommand{\btheta}{\boldsymbol{\theta}}

\newcommand*\R{\mathbb{R}}
\newcommand{\F}{\hat{F}_{\boldsymbol{\theta}}}
\newcommand{\f}{\hat{f}_{\boldsymbol{\theta}}}

\title{\LARGE{\vspace{-.65in}\textbf{Deep learning of dynamics and signal-noise decomposition with time-stepping constraints}}\vspace{-.25in}}
\author{\normalsize{Samuel H. Rudy$^{1*}$, J. Nathan Kutz$^1$, Steven L. Brunton$^2$}\\
\footnotesize{$^1$ Department of Applied Mathematics, University of Washington, Seattle, WA 98195, United States}\\
\footnotesize{$^2$ Department of Mechanical Engineering, University of Washington, Seattle, WA 98195, United States\vspace{-.8in}}
}
\date{}
\begin{document}
\maketitle

\blfootnote{$^*$ Corresponding author (shrudy@uw.edu).\\ \noindent \textbf{Python code:}  https://github.com/snagcliffs/RKNN}
\vspace{-.25in}
\begin{abstract}
A critical challenge in the data-driven modeling of dynamical systems is producing methods robust to measurement error, particularly when data is limited.  Many leading methods either rely on denoising prior to learning or on access to large volumes of data to average over the effect of noise.  We propose a novel paradigm for data-driven modeling that simultaneously learns the dynamics and estimates the measurement noise at each observation.  
By constraining our learning algorithm, our
method explicitly accounts for measurement error in the map between observations, treating both the measurement error and the dynamics as unknowns to be identified, rather than assuming idealized noiseless trajectories.  
We model the unknown vector field using a deep neural network, imposing a Runge-Kutta integrator structure to isolate this vector field, even when the data has a non-uniform timestep, thus constraining and focusing the modeling effort.   
We demonstrate the ability of this framework to form predictive models on a variety of canonical test problems of increasing complexity and show that it is robust to substantial amounts of measurement error.  We also discuss issues with the generalizability of neural network models for dynamical systems and provide open-source code for all examples.\\

\vspace{-0.05in}
\end{abstract}


\vspace{-.25in}
\section{Introduction}

Dynamical systems are ubiquitous across nearly all fields of science and engineering.  
When the governing equations of a dynamical system are known, they allow for forecasting, estimation, control, and the analysis of structural stability and bifurcations.  
Dynamical systems models have historically been derived via first principles, such as conservation laws or the principle of least action, but these derivations may be intractable for complex systems, in cases where mechanisms are not well understood, and/or when measurements are corrupted by noise. 
These complex cases motivate automated methods to develop dynamical systems models from 
data, although nearly all such methods are compromised by limited and/or noisy measurements.   
In this work, we show that by constraining the learning efforts to a time-stepping structure, the underlying dynamical system can be disambiguated from noise, allowing for efficient and uncorrupted model discovery.

Data-driven system identification has a rich history in science and engineering~\cite{Juang1994book,ljung:book,billings2013nonlinear}, and recent advances in computing power and the increasing rate of data collection have led to renewed and expanded interest.  
Early methods identified linear models from input--output data based on the minimal realization theory of Ho and Kalman~\cite{Ho1965aac}, including the eigensystem realization algorithm (ERA)~\cite{Juang1985jgcd, Longman1989jgcd}, and the observer/Kalman filter identification (OKID)~\cite{Juang1991nasatm,Phan1992jas,Phan1993jota}, which are related to the more recent dynamic mode decomposition (DMD)~\cite{schmid:2010,Tu2014jcd,kutz2016dynamic}.  
OKID explicitly accounts for measurement noise by simultaneously identifying a de-noising Kalman filter~\cite{Kalman1960jfe} and the impulse response of the underlying noiseless system, providing considerable noise robustness for linear systems from limited measurements.  

Several approaches have also been considered for learning interpretable nonlinear models.
Sparse regression techniques have been used to identify exact expressions for nonlinear ordinary differential equations \cite{Wang2011prl,Brunton2016,Loiseau2017jfm,schaeffer2017extracting,schaeffer2017sparse,schaeffer2018extracting, tran2017exact,mangan2017model}, partial differential equations \cite{schaeffer2017learning,rudy2017data}, and stochastic differential equations \cite{boninsegna2017sparse}.
Recent work has also used Gaussian process regression to obtain exact governing equations \cite{raissi2017hidden}.
These methods provide interpretable forms for the governing equations, but rely on libraries of candidate functions and therefore have difficulty expressing complex dynamics.  
Symbolic regression \cite{Bongard2007pnas,Schmidt2009science} allows for more expressive functional forms at the expense of increased computational cost; these methods have been used extensively for automated inference of systems with complex dynamics~\cite{Schmidt2011pb,Daniels2015naturecomm,Daniels2015plosone}.  
Several of the aforementioned techniques have specific considerations for measurement noise.  In particular, the sparse regression methods in \cite{Brunton2016,rudy2017data} use smoothing methods in their numerical differentiation schemes, \cite{tran2017exact} identifies highly corrupted measurements, \cite{schaeffer2017sparse} attenuates error using integral terms, and \cite{raissi2017hidden} naturally treats measurement noise by representing data as a Gaussian process.

Although the interpretable nonlinear models above are appealing, there are currently limitations to the flexibility of functional forms and the number of degrees of freedom that can be modeled.  
There has been considerable and growing recent interest in leveraging powerful black-box machine learning techniques, such as deep learning, to model increasingly complex systems.  
Deep learning is particularly appealing because of the ability to represent arbitrarily complex functions, given enough training data of sufficient variety and quality.   
Neural networks have been used to model dynamical systems for decades~\cite{chen1990non, gonzalez1998identification,Milano2002jcp}, although recent advances in computing power, data volumes, and deep learning architectures have dramatically improved their capabilities.  Reccurrent neural networks naturally model sequential processes and have been used for forecasting \cite{vlachas2018data, pan2018long,lu2018attractor, pathak2018hybrid, pathak2017using} and closure models for reduced order models \cite{wan2018data, pan2018data}.  Deep learning approaches have also been used recently to find coordinate transformations that make strongly nonlinear systems approximately linear, related to Koopman operator theory~\cite{Takeishi2017nips,Yeung2017arxiv,Wehmeyer2017arxiv,Mardt2017arxiv,Lusch2017arxiv}.  Feed-forward networks may also be used in conjunction with classical methods in numerical analysis to obtain discrete timesteppers \cite{raissi2018multistep, raissi2017physics1, raissi2017physics2, raissi2018deep, gonzalez1998identification}. 
Many of these modern identification and forecasting methods may be broadly cast as nonlinear autoregressive moving average models with exogenous inputs (NARMAX) \cite{chen1989representations, billings2013nonlinear} with increasingly sophisticated interpolation models for the dynamics. 
NARMAX explicitly accounts for noise and forcing terms in prediction, and estimates exact values of measurement noise by alternating between learning model parameters and noise \cite{billings2013nonlinear}.

The methods presented in \cite{raissi2018multistep, gonzalez1998identification} are highly structured nonlinear autoregressive models.  
In each case, the dynamics are represented using classical methods from numerical analysis where neural networks are used to interpolate the underlying vector field and time-stepping is performed by multistep or Runge-Kutta methods.  This framework is highly compelling because it allows for the interpolation of the vector field as opposed to the discrete map, which is most likely a more complicated function.  However, these models do not explicitly account for measurement noise.  
In this work, we expand on \cite{raissi2018multistep, gonzalez1998identification} to explicitly account for measurement noise by constructing a framework for learning measurement noise in tandem with the dynamics, rather than a sequential or alternating optimization.  Taken together, these innovations provide considerable robustness to measurement noise and reduce the need for vast volumes of data.

\subsection{Contribution of this work}
The principal contribution of this work is to introduce a new paradigm for learning governing equations  from noisy time-series measurements where we account for measurement noise explicitly in a map between successive observations. 
By constraining the learning process inside a numerical timestepping scheme, we can improve the ability of automated methods for model discovery by cleanly separating measurement error from the underlying state while simultaneously learning a neural representation of the governing equations.  
A Runge-Kutta integration scheme is imposed in the optimization problem to focus the neural network to identify the continuous vector field, even when the data has a variable timestep.  
Our method yields predictive models on a selection of dynamical systems models of increasing complexity even with substantial levels of measurement error and limited observations. 
We also highlight the inherent risks in using neural networks to interpolate governing equations.  
In particular, we focus on overfitting and the challenge of fitting dynamics away from a bounded attractor.  Trade-offs between black-box representations, such as neural networks, and sparse regression methods are considered.

The paper is organized as follows:  In Sec.~\ref{methods} we outline the problem formulation, introduce our computational graph for defining a map between successive observations, describe the optimization problem used to estimate noise and dynamics, and discuss various measures of error.  In Sec.~\ref{results} we apply our method to a selection of test problems of increasing complexity.  Section \ref{caution} discusses  pitfalls of using a black box machine learning technique such as a neural network as a model for the underlying vector field.  Several examples are given to demonstrate radical differences between true and learned vector fields, even when test trajectories are accurately predicted.  Section \ref{discussion} provides conclusions and a discussion about further directions.

\begin{figure*}
\vspace{-.1in}
\framebox[\textwidth]{
\noindent\begin{minipage}[h]{.495\textwidth}
\section*{\hspace{2 mm} Nomenclature}
\begin{tabular}{ll}
$t_j$ & Timestep $j$ \\
$\mathbf{x}_j$ & State at time $t_j$ \\
$\boldsymbol{\nu}_j$ & Measurement error at time $t_j$ \\
$\mathbf{y}_j$ & Measurement at time $t_j$ \\

$n$ & Dimension of state space\\
$m$ & Number of timesteps\\

$\mathbf{X}$ & State matrix in $\R^{n \times m}$\\
$\mathbf{N}$ & Measurement error matrix in $\R^{n \times m}$\\
$\mathbf{Y}$ & Measurement matrix in $\R^{n \times m}$\\

$f$ & True underlying vector field\\
$F$ & Discrete-time flow map\\

\end{tabular}
\end{minipage}
\begin{minipage}[h]{.495\textwidth}
\begin{tabular}{ll}

$\btheta$ & Parameterization of neural network $\f$\\
$\f$ & Data-driven approximation of $f$\\
$\F$ & Data-driven approximation of $F$\\
$q$ & \# of forward and backward steps\\

$E_N$ & Measurement noise error\\
$E_{f}$ & Vector field error\\
$E_{F}$ & Forward orbit error\\

$\mathcal{L}$ & Loss function for $\btheta$ and $\mathbf{N}$\\

\end{tabular}

\vspace{1 mm}
Color codes on figures: \\
{\color{blue} Observed data}\\
{\color{red} Learned parameters}\\
{\color{green} Model predicted trajectory}
\end{minipage}
}
\vspace{-0.15in}
\end{figure*}

\section{Methods}\label{methods}
\begin{figure}
\centering
\includegraphics[width=\textwidth]{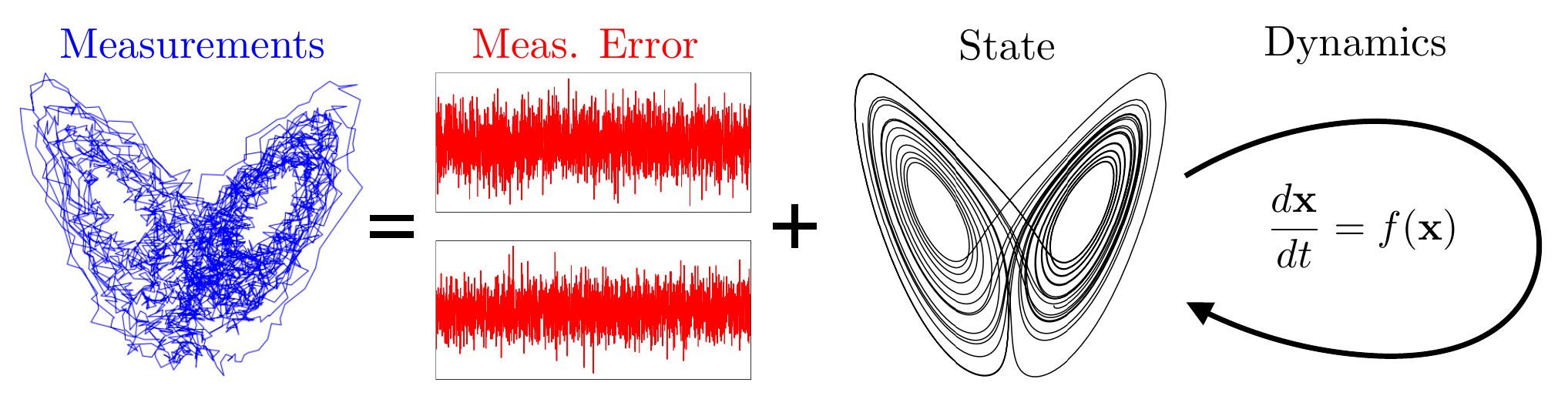}
\caption{Measurements of a dynamical system may be split into measurement error and an underlying state that is governed by a dynamical system.}
\label{fig:basic_schematic}
\end{figure}

The methods presented in this work add to a growing body of literature concerned with the machine learning of dynamics from data.  
A fundamental problem in the field has been formulating methods that are robust to measurement noise.  
Many previous works have either relied on smoothing techniques or large datasets to average out the noise.  
We propose a novel approach by treating measurement error as a latent variable relating observations and a true state governed by dynamics, as in Fig.~\ref{fig:basic_schematic}.  
Our method generates accurate predictive models from relatively small datasets corrupted by large amounts of noise by explicitly considering the measurement noise as part of the model instead of smoothing the data.  
The computational methodology simultaneously learns pointwise estimates of the measurement error, which are subtracted from the measurements to estimate the underlying state, as well as a dynamical model to propagate the true state.  
Section \ref{problem_formulation} provides an overview of the class of problems we consider, Sec. \ref{comp_framework} discusses our computational approach, and Sec. \ref{measuring_error} provides metrics to evaluate the accuracy of the data-driven model.

\subsection{Problem formulation} \label{problem_formulation}

We consider a continuous dynamical system of the form,
\begin{equation}
\frac{d}{dt} \mathbf{x}(t) = f\left(\mathbf{x}(t)\right),
\label{eq:basic_ode}
\end{equation}
where $\mathbf{x} \in \mathbb{R}^n$ and $f$ is a Lipschitz continuous vector field.  This work assumes that $f$ is unknown and that we have a set of measurements $\mathbf{Y} = [\mathbf{y}_1, \hdots, \mathbf{y}_m]$ with $\mathbf{y}_j \in \mathbb{R}^n$ representing the true state at time $t_j$ corrupted by measurement noise $\boldsymbol{\nu}_j$:
\begin{equation}
\mathbf{y}_j = \mathbf{x}_j + \boldsymbol{\nu}_j.
\label{eq:noise_assumtion}
\end{equation}
The discrete-time map from $\mathbf{x}_j$ to $\mathbf{x}_{j+1}$ may be written explicitly as
\begin{equation}
\mathbf{x}_{j+1} = F(\mathbf{x}_j) = \mathbf{x}_j + \int_{t_j}^{t_{j+1}} f(\mathbf{x}(\tau))\,d\tau .
\label{eq:basic_timestepper}
\end{equation}
Many leading methods to approximate the dynamics (e.g., $f$ or $F$) from data involve a two-step procedure, where pre-processing methods are first applied to clean the data by filtering or smoothing noise, followed by the second step of fitting the dynamics.  
In contrast, we seek to simultaneously approximate the function $f$ and obtain estimates of the measurement error $\boldsymbol{\nu}_j$ at each timestep.  
Specifically, we consider both the dynamics and the measurement noise as unknown components in a map between successive observations:
\begin{equation}
\mathbf{y}_{j+i} = \overbrace{F^i(\underbrace{\mathbf{y}_j - \boldsymbol{\nu}_j}_{\mathbf{x}_j})}^{\mathbf{x}_{j+i}} + \boldsymbol{\nu}_{j+i}.
\label{eq:obs_to_obs_true}
\end{equation}
Starting with the observation $\mathbf{y}_j$, if the measurement noise $\boldsymbol{\nu}_j$ is known, it may be subtracted to obtain the true state $\mathbf{x}_j$.  The flow map is applied $i$ times to the state $\mathbf{x}_j$ to obtain the true state at timestep $\mathbf{x}_{j+i}$, and adding the measurement noise back will yield the observation $\mathbf{y}_{j+i}$.  Given that $\mathbf{y}_j$ is observed, there are two unknown quantities in \eqref{eq:obs_to_obs_true}: the dynamics $F$ and the measurement noise $\mathbf{N} = [\boldsymbol{\nu}_1,\hdots , \boldsymbol{\nu}_m]$.  
We leverage the fact that \eqref{eq:obs_to_obs_true} must hold for all pairs of observations including $i<0$ and enforce consistency in our estimate of $\boldsymbol{\nu}_j$ to separate dynamics from measurement noise.
In this framework, consistency of governing dynamics helps decompose the dataset into state and measurement error, while explicit estimates of measurement error simultaneously allow for a more accurate model of the dynamics. 
In addition, we approximate $F$ with a Runge-Kutta scheme and focus on modeling the continuous dynamics $f$, which is generally simpler than $F$.

\subsection{Computational framework}\label{comp_framework}

\begin{figure}
\centering
\includegraphics[width=\textwidth]{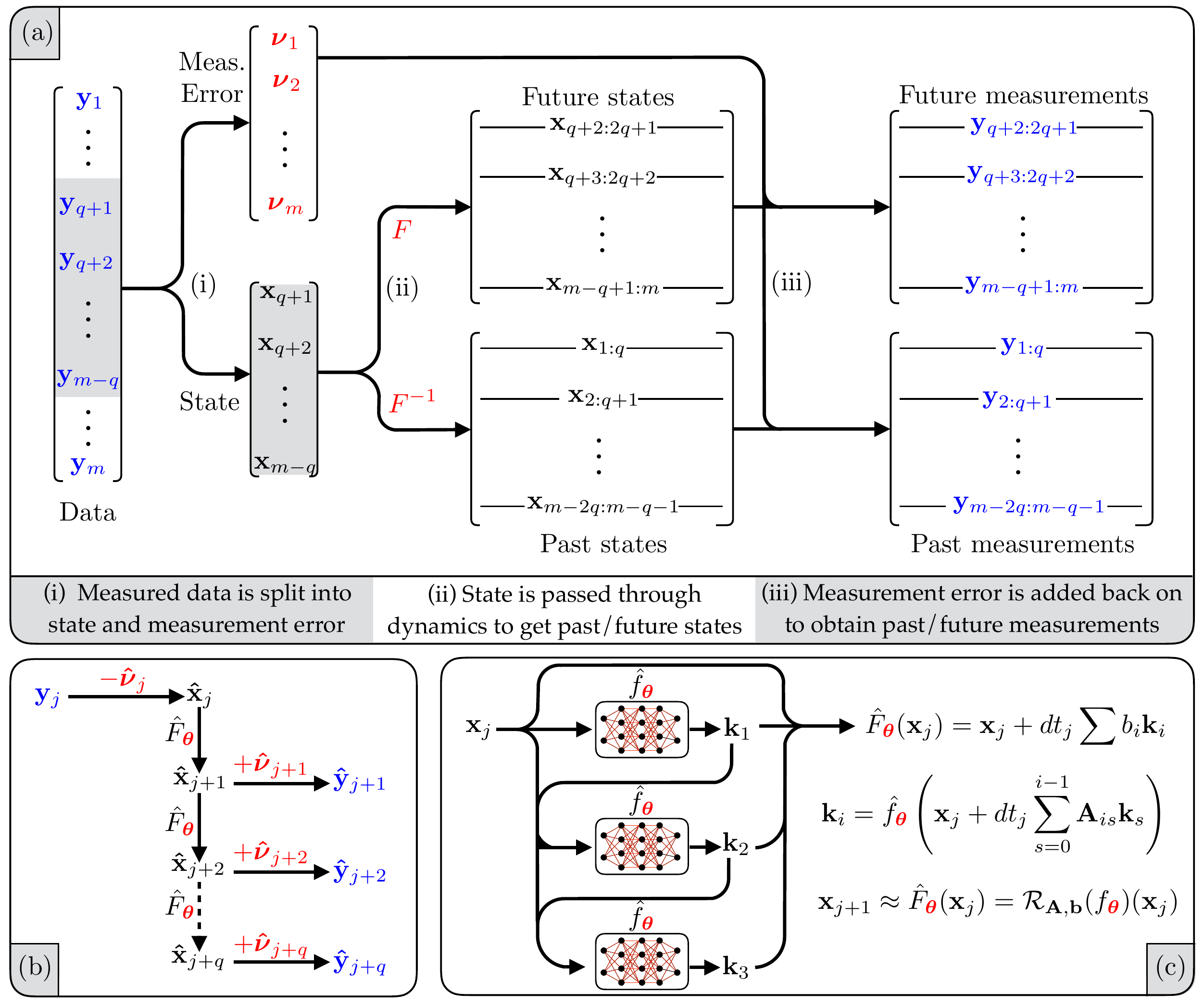}
\caption{ Schematic of de-noising framework for constrained learning of nonlinear system identification.  {\bf(a)} Flow of full dataset through \eqref{eq:obs_to_obs_true} to find forward and backward measurements.  Measurements are decomposed into learned measurement error and state, the state is propagated forward and backward using learned dynamics, and measurement error is added back to find forward and backward measurements and evaluate model accuracy. {\bf(b)} Example of data-driven estimate of forward propagation of single observation from graph in panel a. {\bf(c)} Example 3-step Runge-Kutta scheme to approximate the map $F$ and isolate the vector field, which is modeled by a feed-forward neural network.}
\label{fig:main_fig}
\end{figure}

Here we outline a method for estimating both the unknown measurement noise and model in \eqref{eq:obs_to_obs_true}, shown in Fig.~\ref{fig:main_fig}.  
We let $\mathbf{\hat{N}} = [\boldsymbol{\hat{\nu}}_1,\hdots , \boldsymbol{\hat{\nu}}_m] \in \R^{n \times m}$ and $\F : \R^n \to \R^n$ denote the data-driven approximation of the measurement noise and discrete-time dynamics, respectively.  
We develop a cost function to evaluate the fidelity of model predictions in comparison with the observed data and then we solve an optimization problem to refine estimates for both the measurement noise $\mathbf{\hat{N}}$ and the parameterization of $\F$.  
The map $\F$ is modeled by a Runge-Kutta time-stepper scheme that depends on a continuous vector field, $f$, which is modeled by a feed-forward network.  
 
We construct a data-driven model for the vector field $f$ using a feed-forward neural network:
\begin{equation}
\f(\mathbf{x}) = \left(\prod_{i=1}^l C_{g_i} \right) (\mathbf{x})\,\, \text{ where }\,\, g_i(\mathbf{x}) = \sigma_i(\mathbf{W}_i\mathbf{x}+\mathbf{c}_i)\, ,
\label{eq:feed_forward_nn}
\end{equation}
where $C_g$ is taken to be the composition operator with $g$, each $\sigma_i$ is a nonlinear activation function, and $\btheta = \{\mathbf{W}_i,\mathbf{c}_i\}_{i=1}^l$ denotes the neural network's parameterization.  A large collection of activation functions exist in the literature, with the standard choice for deep networks being the rectified linear unit (ReLU) \cite{goodfellowdeep}.  For interpolating the vector field, we prefer a smooth activation function and instead use the exponential linear unit, given by 
\begin{equation}
\sigma(\mathbf{x}) = \left\{ \begin{aligned}
&e^\mathbf{x}-1, \text{ for } \mathbf{x}\leq 0 \\
&\mathbf{x}, \hspace{9 mm} \text{for } \mathbf{x}> 0
\end{aligned}\right. 
\label{eq:elu}
\end{equation}
evaluated component-wise for $\sigma_1, \hdots, \sigma_{l-1}$, and the identity for $\sigma_l$.  Neural networks are typically trained by passing in known pairs of inputs and outputs, but doing so for $\f$ would require robust approximations of the time derivative from noisy time series, which can be prohibitively difficult.  

To avoid numerical differentiation of noisy data, we embed $\f$ into a time-stepper to obtain a discrete map approximating the flow of the dynamical system, \eqref{eq:basic_timestepper}.  We use explicit Runge-Kutta schemes because they allow us to make predictions from a single state measurement and are easily formed as an extension of the neural network for $\f$, shown in Fig.~\ref{fig:main_fig} (c).  Runge-Kutta schemes for autonomous systems are uniquely defined by weights $\mathbf{A} \in \R^{p\times p}, \mathbf{b} \in \R^p$, where $p$ the the number of intermediate steps \cite{leveque2007finite}.  Given $\mathbf{A},\mathbf{b}$, we let $\mathcal{R}_{\mathbf{A},\mathbf{b}}$ denote the operator induced by a time-stepper with parameters $\mathbf{A},\mathbf{b}$ mapping a vector field to a discrete flow map.  Applying $\mathcal{R}_{\mathbf{A},\mathbf{b}}$ to an approximation of the vector field $\f$ gives an approximate flow map,
\begin{equation}
\mathbf{\hat{x}}_{j+1} = \F (\mathbf{x}_j) = \mathcal{R}_{\mathbf{A},\mathbf{b}}(\f)(\mathbf{x}_j) \approx F(\mathbf{x}_j).
\label{eq:flow_map}
\end{equation}
Incorporating estimates for the measurement error at each timestep extends \eqref{eq:flow_map} to a data-driven map from the observation $\mathbf{y}_j$ to $\mathbf{y}_{j+i}$, approximating the exact map in \eqref{eq:obs_to_obs_true}:
\begin{equation}
\mathbf{\hat{y}}_{j+i} = \F^i(\mathbf{y}_j - \boldsymbol{\hat{\nu}}_j) + \boldsymbol{\hat{\nu}}_{j+i}.
\label{eq:obs_map}
\end{equation}
Discrepancies between $\mathbf{\hat{y}}_{j+i}$ and $\mathbf{y}_{j+i}$ will result from inaccurate estimates of $\boldsymbol{\nu}$, $f$ and from numerical error in the timestepping scheme.  For a single pair of observations $\mathbf{y}_j$ and $\mathbf{y}_{j+i}$ the squared $L^2$ difference given the parameterization $\btheta$ and a noise estimate $\mathbf{\hat{N}}$ is,
\begin{equation}
\mathcal{L}_{j,i}(\btheta, \mathbf{\hat{N}}, \mathbf{Y}) = \left\| \mathbf{y}_{j+i} -\left( \F^i(\mathbf{y}_j - \boldsymbol{\hat{\nu}}_j) + \boldsymbol{\hat{\nu}}_{j+i} \right)  \right\|_2^2.
\label{eq:single_step_loss}
\end{equation}
Summing the error in \eqref{eq:single_step_loss} over all pairs of observations will result in a computationally stiff and intractable optimization problem.   Moreover, chaotic dynamics will render observations statistically uncorrelated as they separate in time.  Instead we formulate a global evaluation metric by summing \eqref{eq:single_step_loss} over pairs $j, j+i$ in a local neighborhood with $|i| \leq q$.  
We also weight each pair with exponentially decreasing importance in $|i|$, given by $\omega_i$, according to the expected accumulation of error given inaccuracies in $\boldsymbol{\hat{\nu}}_j$.  A careful discussion of the weighting is given in Appendix A.  
The resulting cost function is,
\begin{equation}
\mathcal{L} (\btheta, \mathbf{\hat{N}}, \mathbf{Y}) = \sum_{j={q+1}}^{m-q} \sum_{i=-q}^q \omega_i \left\| \mathbf{y}_{j+i}- \left( \F^i(\mathbf{y}_j - \boldsymbol{\hat{\nu}}_j) + \boldsymbol{\hat{\nu}}_{j+i} \right)  \right\|_2^2.
\label{eq:cost_no_regularization}
\end{equation}

The cost function in \eqref{eq:cost_no_regularization} has a global minimum with the trivial solution $\f = 0$ and estimated measurement error $\boldsymbol{\nu}_j = \mathbf{y}_j - \mathbf{y}_1$ accounting for all observed dynamics.  In practice, locally minimizing solutions to \eqref{eq:cost_no_regularization} are nontrivial but do not result in accurate models.  Penalizing the magnitude of $\mathbf{\hat{N}}$ as well as the weights in the neural network results in a much more robust loss function:
\begin{equation}
\mathcal{L} (\btheta, \mathbf{\hat{N}}, \mathbf{Y}) = \sum_{j={q+1}}^{m-q} \sum_{i=-q}^q \omega_i \left\| \F^i(\mathbf{y}_j - \boldsymbol{\hat{\nu}}_j) + \boldsymbol{\hat{\nu}}_{j+i} - \mathbf{y}_{j+i} \right\|_2^2 + \gamma \|\mathbf{\hat{N}}\|_F^2 +\beta \sum_{i=1}^l \|\mathbf{W}_i\|_F^2.
\label{eq:full_cost}
\end{equation}
The regularization term for $\mathbf{\hat{N}}$ in \eqref{eq:full_cost} makes the trivial solution $\f = 0$ highly costly and encourages the neural network to fit dynamics to a time series close to the observations.  Penalizing the weights of the neural network discourages overfitting and is particularly important for larger networks.  
Neural network parameters, as well as explicit estimates of the measurement error, are obtained by minimizing \eqref{eq:full_cost} using the quasi-Newton method L-BFGS-B \cite{zhu1997algorithm} implemented in SciPy.

\subsection{Measuring error}\label{measuring_error}

In addition to the loss function derived in Sec.~\ref{comp_framework}, we use several other metrics to evaluate the accuracy of models produced by minimizing \eqref{eq:full_cost}.  
It is possible to evaluate these metrics for the problems considered in this work because the equations and measurement noise are both known.  
Although these are not generally available for a new dataset, they provide a quantitative basis for comparing performance with different quantities of noise, volumes of data, and timestepping schemes.  
To evaluate the accuracy of $\f$, we use the relative squared $L^2$ error between the true and data-driven vector fields,
\begin{equation}
E_{f}(\f) \, = \dfrac{\sum_{j=1}^m \|f(\mathbf{x}_j) - \f(\mathbf{x}_j)\|_2^2}{\sum_{j=1}^m \|f(\mathbf{x}_j)\|_2^2}.
\label{eq:trajectory_vec_field_error}
\end{equation}
Here the error is only evaluated along the noiseless training data.  
This results in substantially lower error than if the two vector fields were compared on a larger set away from the training data, as will be discussed in Sec.~\ref{caution}.

The other learned quantity is an estimate for the measurement noise $\boldsymbol{\nu}_j$, or equivalently the de-noised state $\mathbf{x}_j$, at each timestep.  The mean $L^2$ difference between the true and learned measurement error is,
\begin{equation}
E_{N}(\mathbf{\hat{N}}) \, = \, \frac{1}{m} \sum_{j=1}^m \|\boldsymbol{\nu}_j - \boldsymbol{\hat{\nu}}_j \|_2^2 = \, \frac{1}{m}\sum_{j=1}^m \|\mathbf{x}_j - \mathbf{\hat{x}}_j \|_2^2.
\label{eq:noise_error}
\end{equation}

The $L^2$ distance between the true state of the training data and the forward orbit of $\mathbf{x}_1$ as predicted by $\F$ is computed as,
\begin{equation}
E_{F}(\F,\mathbf{x}_1) = \frac{1}{\|\mathbf{X}\|_F^2} \sum_{j=1}^{m-1} \left\| \mathbf{x}_j - \F^j (\mathbf{x}_1) \right\|_2^2 .
\label{eq:reconstruction_error}
\end{equation}
The last error \eqref{eq:reconstruction_error} is the most stringent and may be highly sensitive to small changes in the dynamics that reflect numerical error more than inaccuracies.  
$E_{F}$ will not yield informative results for dynamics evolving on a chaotic attractor or slowly diverging from an unstable manifold.  
We therefore only consider it on the first example of a damped oscillator.

\section{Results}\label{results}
In this section we test the performance of the methods discussed in Sec.~\ref{methods} on a range of canonical problems of increasing complexity.  
To demonstrate robustness, we consider each problem with varying levels of corruption, meant to replicate the effects of measurement noise.  
In most cases, independent and identically distributed Gaussian noise is added to each component of the dataset with zero mean and amplitude equal to a given percent of the standard deviation of the data: 
\begin{equation}
\boldsymbol{\nu} \sim \mathcal{N}(0, \mathbf{\Sigma}_\mathbf{N}^2),\, \hspace{1 cm}\mathbf{\Sigma}_\mathbf{N} = \dfrac{\text{noise \%}}{100}\text{diag}(\mathbf{\Sigma}_\mathbf{Y}) .
\label{eq:artificial_noise}
\end{equation}
For the Lorenz equation, the measurement noise is drawn from a Student's T distribution.  
In each case, an initial estimate for noise is obtained by a simple smoothing operation performed on the raw data.  
While not necessary, this was found to speed up the optimization routine.  
Weights for the neural networks are initialized using the Xavier initialization native to TensorFlow \cite{glorot2010understanding}.

\begin{figure}
\centering
\includegraphics[width=\textwidth]{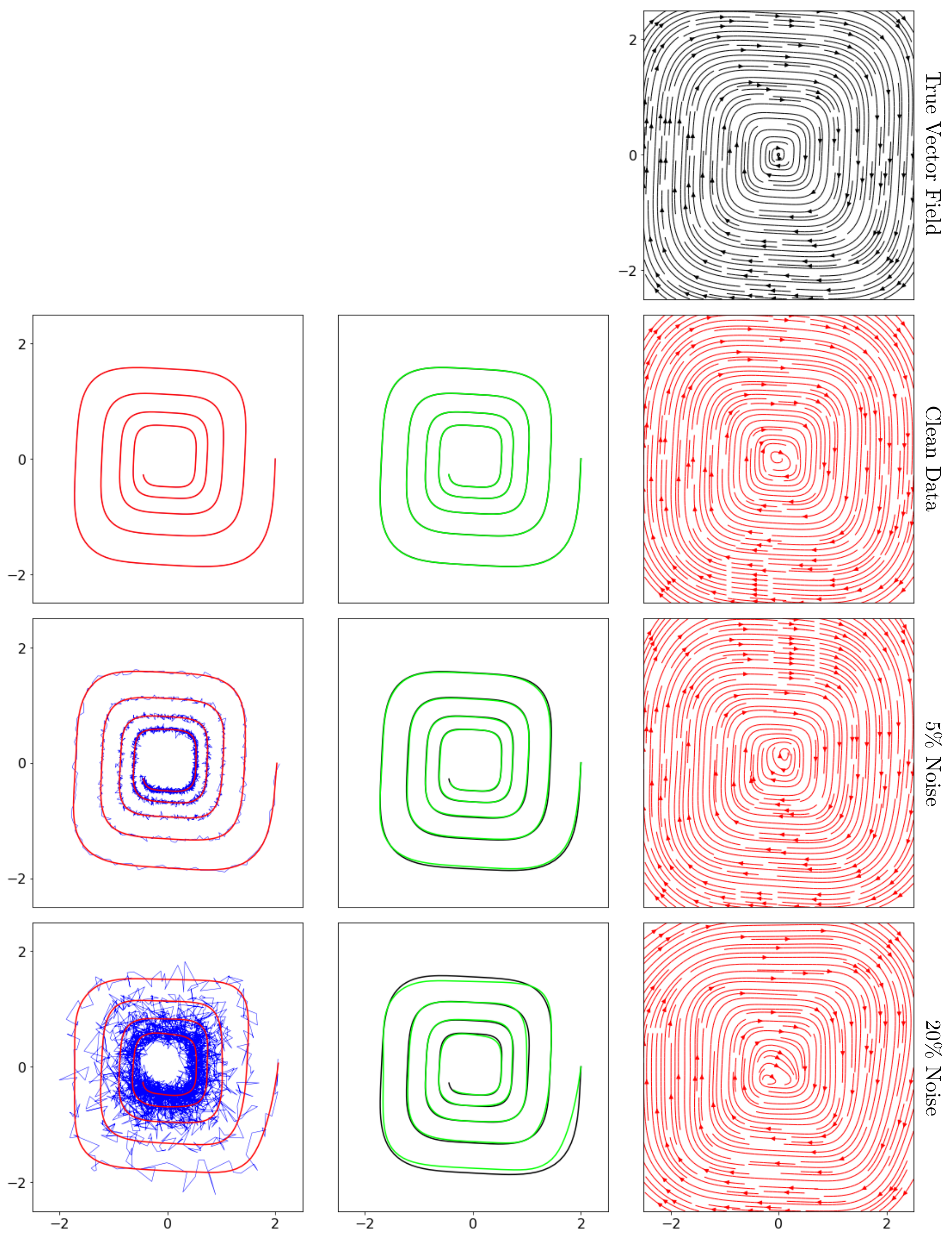}
\begin{picture}(0,0)
\put(-220,560){\parbox{.55\linewidth}{
\noindent \caption{Results for cubic oscillator with increasing magnitudes of measurement noise.  Left column: Observations $\mathbf{Y}$ in blue and learned state $\mathbf{Y}-\mathbf{\hat{N}}$ in red.  Middle column: True state $\mathbf{X}$ in black and forward orbit of $\mathbf{x}_1$ according to data driven dynamics $\F$ in green.  Right column: Learned neural network representation of vector field $\f$.  True vector field is shown in top right for comparrison.} \label{fig:cubic_osc_main}}}
\end{picture}
\end{figure}
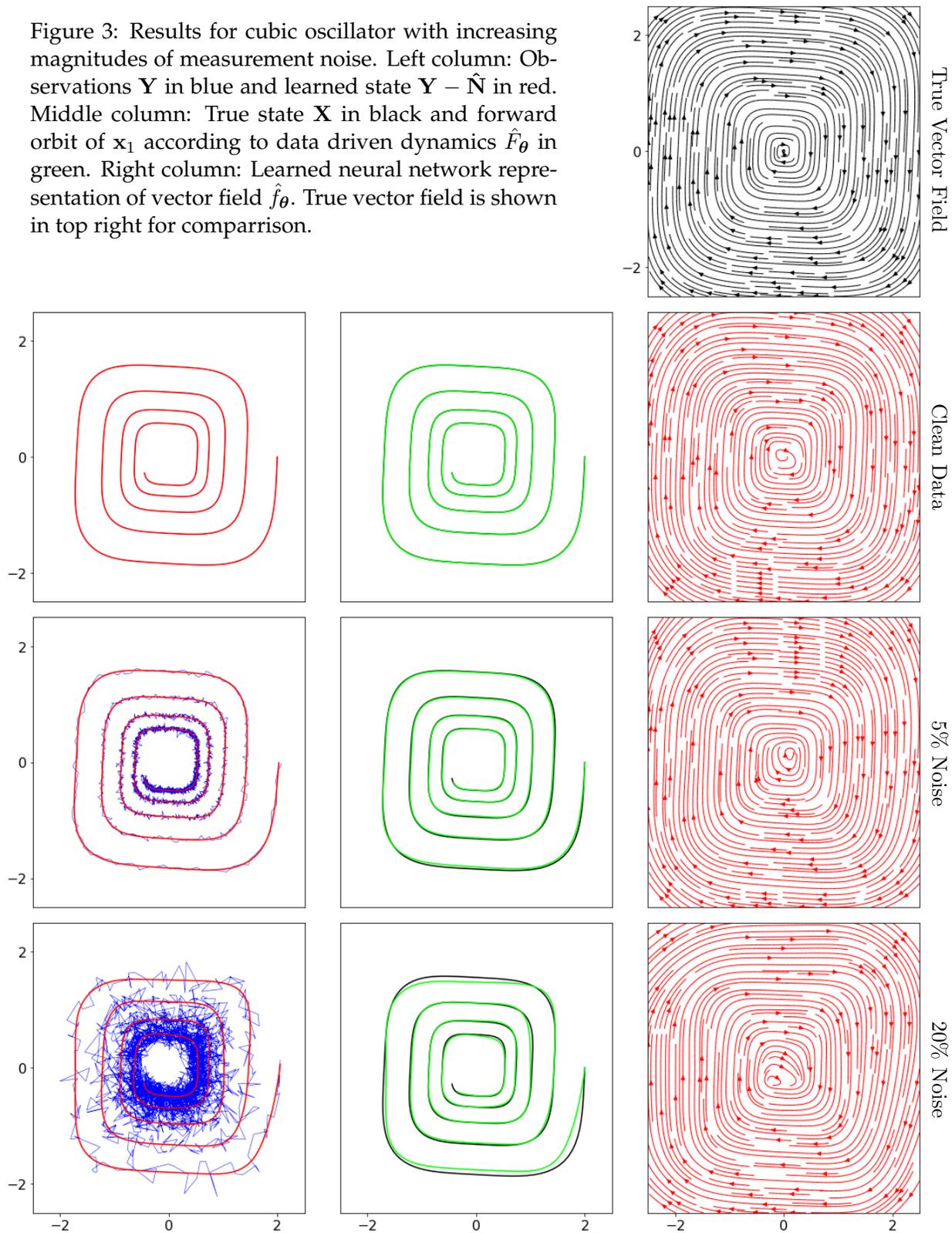
\subsection{Cubic oscillator}
For the first example, we consider the damped cubic oscillator, which has been used as a test problem for several data-driven methods to learn dynamical systems~\cite{Brunton2016,raissi2018multistep}:
\begin{equation}
\begin{aligned}
\dot{x} &= -0.1 x^3 + 2y^3\\
\dot{y} &= -2x^3 - 0.1y^3. 
\end{aligned} \label{eq:cubic_osc}
\end{equation}
We generate $2,500$ snapshots from $t=0$ to $t=25$ via high-fidelity simulation and then corrupt this data with varying levels of artificial noise, ranging from $0$ to $20$ percent. 
Models are trained by embedding a neural network with three hidden layers, containing $32$ nodes each, in a four-step Runge-Kutta scheme.  
For each dataset, we obtain explicit estimates of the measurement noise and a neural network approximation of the vector field in \eqref{eq:cubic_osc}, which is used to integrate a trajectory from the same initial condition to compare the relative error.  
Table \ref{tab:cubic_osc_error} provides a summary of the error metrics from Sec.~\ref{measuring_error} evaluated across various noise levels.  
At higher noise levels, there is a substantial increase in $E_{\F}$ due to a phase shift in the reconstructed solution.
We also tested the method on a dataset with random timesteps, drawn from an exponential distribution, $t_{j+1}-t_j \sim \exp(0.01)$.  
Error in approximating the vector field and noise were similar to the case with a constant timestep, while $E_F$ was significantly lower.  
This suggests future work to perform a careful comparison between the cases of constant and variable timsteps.

\begin{table}
  \begin{center}
    \caption{Error for cubic oscillator model with varying noise.}
    \label{tab:cubic_osc_error}
    \begin{tabular}{|l|c|c|c|c|c|c|}
      \hline
      \% Noise & $0$ & $1$ & $5$ & $10$ & $10$, $dt\sim exp$ & $25$ \\
      \hline
      $E_N$ & $3.077e-7$ & $2.570e-6$ & $3.561e-5$ & $1.292e-4$ & $1.297e-4$ & $7.847e-4$ \\
      $E_{f}$ & $5.187e-5$ & $1.516e-4$ & $7.692e-4$ & $2.436e-3$ & $2.724e-3$ & $1.319e-2$ \\
      $E_{F}$ & $1.150e-4$ & $1.955e-3$ & $3.735e-2$ & $0.2170$ & $5.124e-3$ & $0.6241$ \\
      \hline
    \end{tabular}
  \end{center}
\end{table}

Figure \ref{fig:cubic_osc_main} shows the model predictions and vector fields for increasing amounts of measurement noise.  
In the left column, observations $\mathbf{Y}$ are plotted against the inferred state, $\mathbf{\hat{X}} = \mathbf{Y} - \mathbf{\hat{N}}$.  
The middle column shows the noiseless trajectory $\mathbf{X}$ alongside the forward orbit of $\mathbf{x}_1$ according to the learned timestepper $\F$.  
The learned vector field $\f$ for each noise level is plotted in the right column with the true vector field for reference.  
Results show that the proposed method is highly robust to significant measurement noise.

\begin{figure}
\centering
\includegraphics[width=\textwidth]{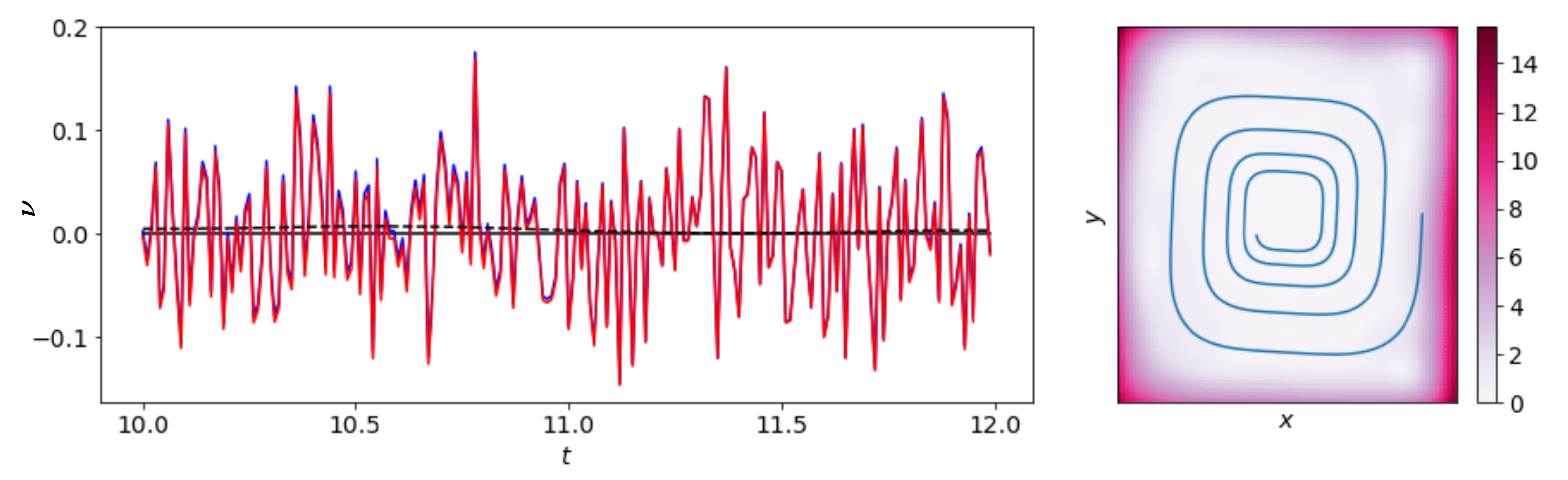}
\vspace{-.3in}
\caption{Results from a single trial of the cubic oscillator with $10\%$ noise.  Left:  True measurement noise (blue), learned measurement noise (red), and error (dotted black line).  Right:  Heat map of vector field approximation error with the noiseless training trajectory plotted for reference.  Note the significant increase in error magnitude in regions far from the training data.}
\label{fig:cubic_osc_single_run}
\end{figure}

Figure \ref{fig:cubic_osc_single_run} shows the error in the approximation of the measurement noise and the vector field for a single time series corrupted by 10\% Gaussian measurement noise.  
The exact measurement noise in the $x$-coordinate is shown alongside the learned measurement noise for $200$ timesteps.  
The $L^2$ error in the approximation of the vector field is also shown with the uncorrupted training trajectory for a spatial reference.  
The vector field error is generally small near the training data and grows considerably near the edges of the domain.

\begin{figure}
\centering
\includegraphics[width=\textwidth]{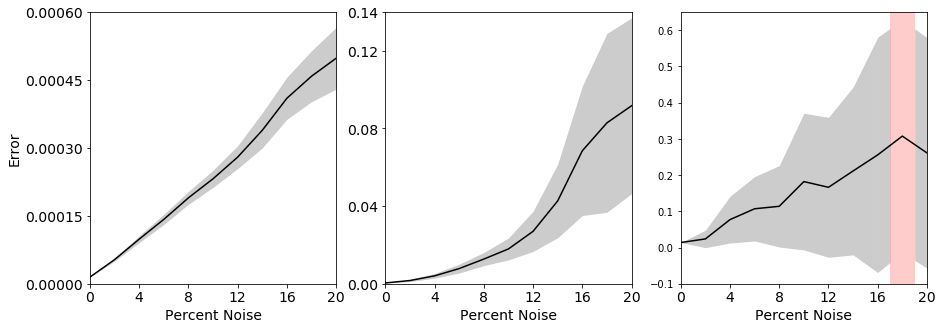}
\vspace{-.3in}
\caption{Average error over 50 trials and standard deviations for the cubic oscillator at varying levels of measurement noise.  Left: Mean $E_N$ and shaded region indicating one standard deviation.  Center:  Mean $E_f$ and standard error.  Right:   Median normalized $L^2$ difference between true trajectory and forward orbit under $\f$, given by $E_F$.  The distribution for $18\%$ noise is omitted due to a single unstable trajectory resulting in non-numeric data; however, the median is reported.}
\label{fig:cubic_osc_errors}
\end{figure}

Exact measures of error will vary depending on the particular instantiation of measurement noise and initialization of model parameters.  
Figure \ref{fig:cubic_osc_errors} shows the average error across fifty trials for each of the metrics discussed in Sec.~\ref{measuring_error}.  
We compute averages for $E_F$ using median rather than mean, since a single trajectory out of the fifty trials with $18\%$ noise was divergent, resulting in non-numerical values.

\subsection{Lorenz system}
\begin{figure}
\centering
\includegraphics[width=.925\textwidth]{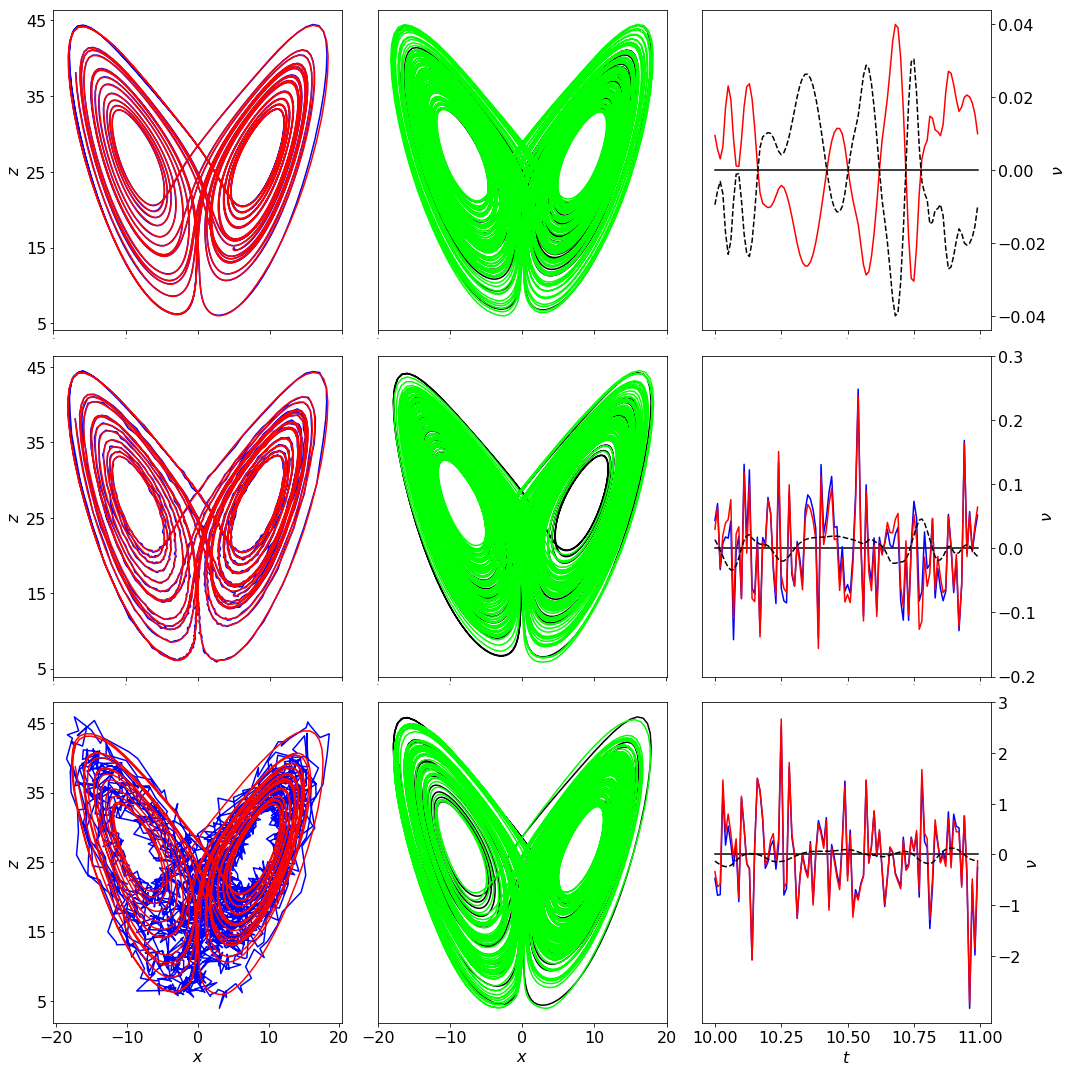}
\vspace{-.15in}
\caption{Results for Lorenz system with increasing magnitudes of measurement noise.  Top row: Clean data.  Second row: 1\% Gaussian noise.  Bottom row: 10\% noise from Student's T distribution with 10 degrees of freedom.  Left column: Observations $\mathbf{Y}$ in blue and learned state $\mathbf{Y}-\mathbf{\hat{N}}$ in red.  Middle column: True state $\mathbf{X}$ in black and forward orbit of $\mathbf{x}_1$ under $\F$ in green.  The prediction is extended to $T_{max}$ five times that of training data.  Right column: True measurement noise $\mathbf{N}$ in blue, learned measurement noise $\mathbf{\hat{N}}$ in red, and error in noise estimate in dashed black.}
\label{fig:lorenz_reconstruction}
\end{figure}

The next example is the Lorenz system, which originated as a simple model of atmospheric convection and became the canonical example for demonstrating chaotic behavior.  
We consider the Lorenz system with the standard parameter set $\sigma = 10$, $\rho = 28$, and $\beta = 8/3$:
\begin{equation}
\begin{aligned}
\dot{x} &= \sigma (y-x)\\
\dot{y} &= x(\rho - z) - y\\
\dot{z} &= xy - \beta z.
\end{aligned} \label{eq:Lorenz}
\end{equation}
The training dataset consists of a single trajectory with $2,500$ timesteps from $t=0$ to $t=25$ with initial condition $(x_0,y_0,z_0) = (5,5,25)$ starting near the attractor.  
The vector field $f$ in \eqref{eq:Lorenz} is modeled by a neural network with three hidden layers containing $64$ nodes each, embedded in a four-step Runge Kutta scheme to approximate $F$. 
Results for several levels of measurement corruption, including noise drawn from a Student's T distribution, are shown in Fig.~\ref{fig:lorenz_reconstruction}.  
Approximation errors for the Lorenz system at varying levels of Gaussian distributed noise are summarized in table \ref{tab:lorenz_error}.

\begin{table}[b]
  \begin{center}
    \caption{Error for Lorenz system with varying noise.}
    \label{tab:lorenz_error}
    \begin{tabular}{|l|c|c|c|c|c|}
      \hline
      \% Noise & $0$ & $1$ & $5$ & $10$ & $15$ \\
      \hline
      $E_N$ & $4.722e-3$ & $4.783e-3$ & $2.670e-2$ & $7.356e-2$ & $0.2248$ \\
      $E_{f}$ & $9.892e-4$ & $9.361e-4$ & $1.988e-3$ & $2.299e-3$ & $6.340e-3$ \\
      \hline
    \end{tabular}
  \end{center}
\end{table}

\begin{figure}
\centering
\includegraphics[width=.9\textwidth]{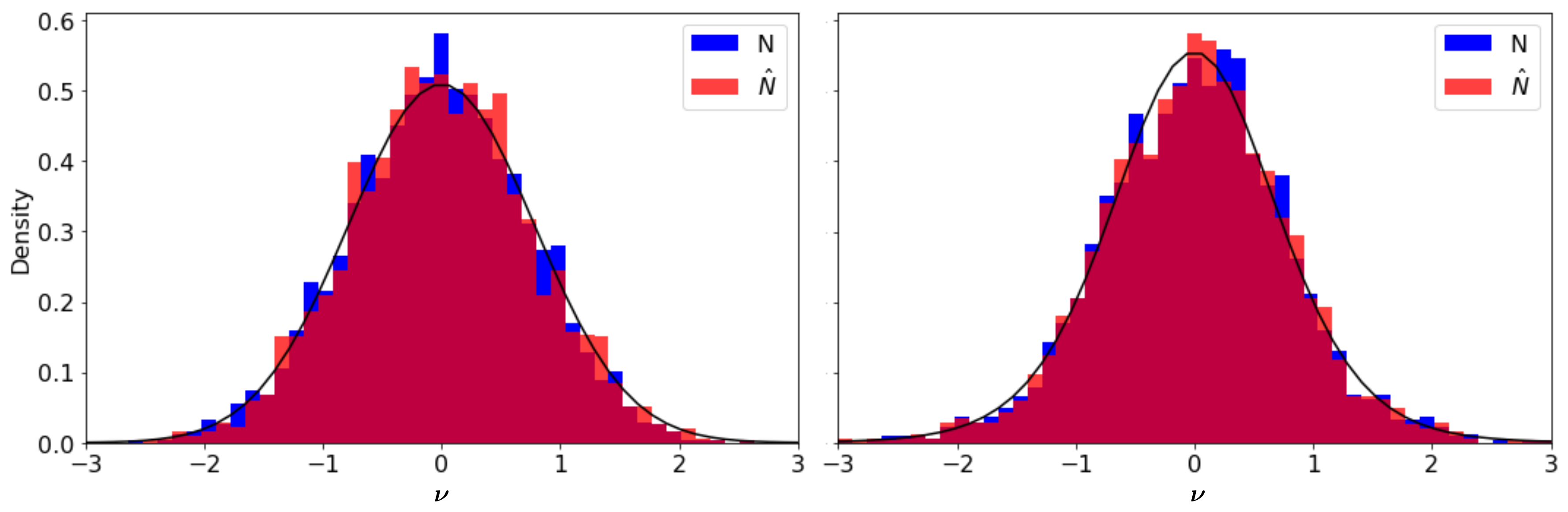}
\caption{Histograms showing true and learned sample distribution of measurement noise with distribution of measurement noise plotted in black.  Left:  Learned noise from Lorenz system with 10\% Gaussian distributed noise.  Right:  Learned noise from Lorenz system with 10\% noise from Student's T distribution with 10 degrees of freedom.}
\label{fig:lorenz_noise_histogram}
\end{figure}

In many cases, it may be important to estimate the distribution of the measurement noise, in addition to the point-wise estimates.  
Figure \ref{fig:lorenz_noise_histogram} shows the true empirical measurement error distribution for the training data along with the distribution of the learned measurement error for the Lorenz system corrupted with either 10\% Gaussian noise or 10\% noise from a Student's T distribution with 10 degrees of freedom; the analytic distribution is also shown for reference. 
In both cases, the approximated error distribution faithfully captures the true underlying distribution of the measurement error.  
Mean, variance, skew, and excess kurtosis of the analytic, empirical, and learned distribution of measurement noise for the x-coordinate are shown in table \ref{tab:lorenz_moments}.

\begin{table}
  \begin{center}
    \caption{Moments of analytic, empirical, and learned measurement noise in $x$-coordinate.}
    \label{tab:lorenz_moments}
    \begin{tabular}{|l|c|c|c|c|c|c|c|c|}
      \hline
      \multicolumn{4}{|c|}{Gaussian Meas. Error} & & \multicolumn{4}{|c|}{Student's T Meas. Error}\\
      \hline
      $\mu$ & $\sigma^2$ & $\gamma_1$ & $\kappa$ & & $\mu$ & $\sigma^2$ & $\gamma_1$ & $\kappa$\\
      \hline
      0 & 1.111 & 0 & 0 & & 0 & 0.6143 & 0 & 1 \\ 
      -0.0145 & 0.5852 & 0.0171 & -0.0920 & & -0.0192 & 0.6143 & 0.0150 & 0.6242\\ 
      -0.0006 & 0.5794 & 0.0093 & -0.0754 & & -0.0003 & 0.6055 & -0.0633 & 0.7718\\
      \hline
    \end{tabular}
  \end{center}
\end{table}

\subsection{Low Reynolds number fluid flow past a cylinder}
As a more complex example, we consider the high-dimensional data generated from a simulation of fluid flow past a circular cylinder at a Reynolds number of $100$ based on cylinder diameter.  
Flow around a cylinder has been a canonical problem in fluid dynamics for decades.  
One particularly interesting feature of the flow is the presence of a Hopf bifurcation occuring at $Re=47$, where the flow transitions from a steady configuration to laminar vortex shedding.  
The low-order modeling of this flow has a rich history, culminating in the celebrated mean-field model of Noack et al.~\cite{Noack2003jfm}, which used Galerkin projection and a separation of timescales to approximate the cubic Hopf nonlinearity with quadratic nonlinearities arising in the Navier-Stokes equations.  
This flow configuration has since been used to test nonlinear system identification~\cite{Brunton2016,Loiseau2017jfm}, and it was recently shown that accurate nonlinear models could be identified directly from lift and drag measurements on the cylinder~\cite{Loiseau2018jfm}.  

\begin{figure}
\centering
\includegraphics[width=0.7\textwidth]{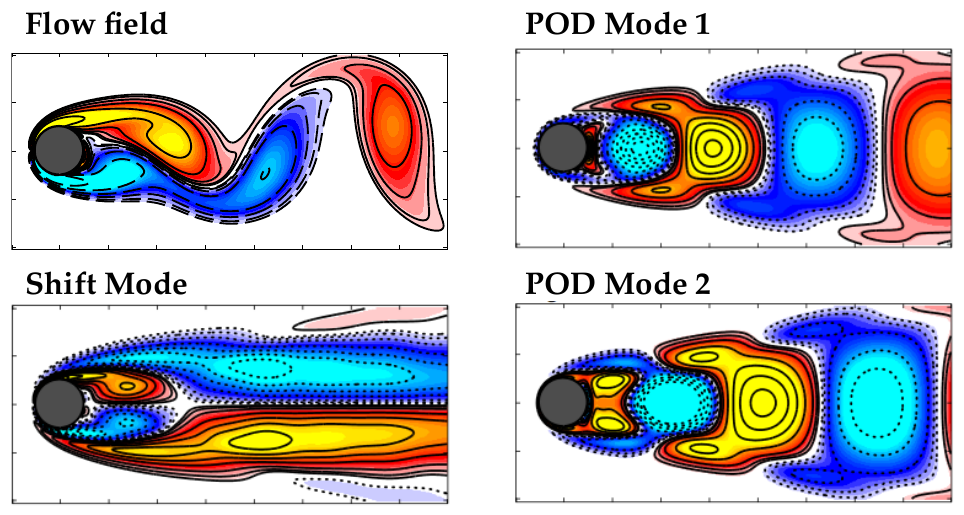}
\caption{Modes from simulation of Navier Stokes equation.  The $x,y,z$ coordinates used in this test problem are time series obtained by projecting the full flow field onto POD mode 1, POD mode 2, and the shift mode, respectively.}
\label{fig:NS_modes}
\end{figure}

\begin{figure}
\centering
\includegraphics[width=\textwidth]{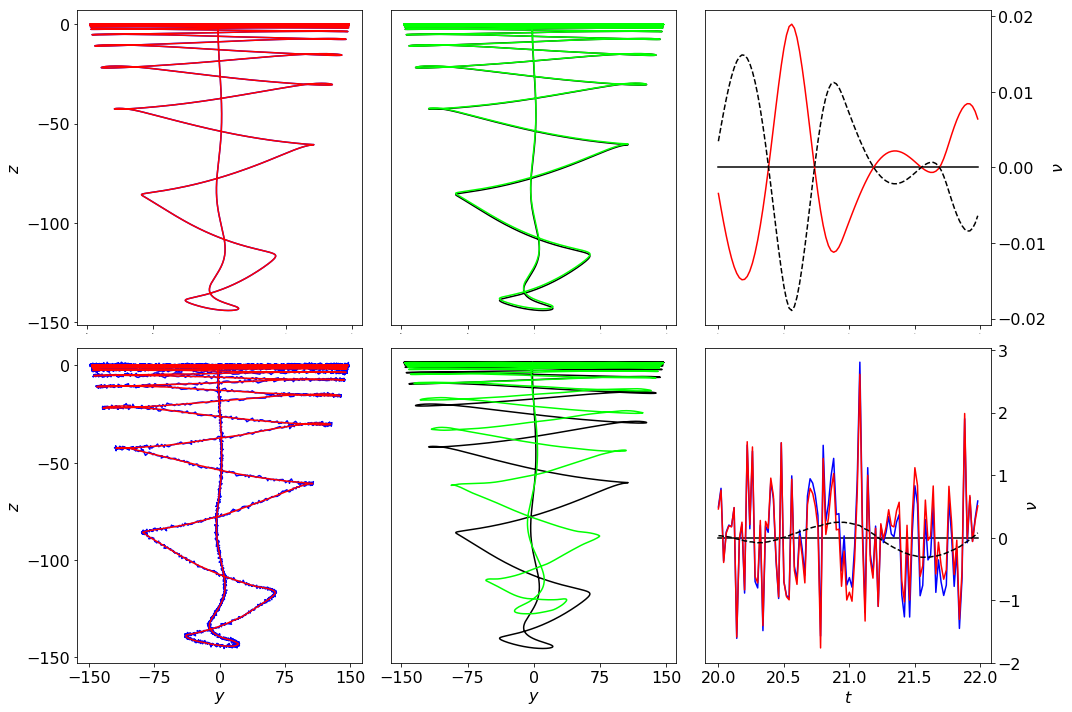}
\caption{Results for reduced basis of flow around a cylinder with increasing magnitudes of measurement noise.  Top row: Clean data.  Bottom row: 1\% Gaussian noise.  Left column: Observations $Y$ in blue and learned state $\mathbf{Y}-\mathbf{\hat{N}}$ in red.  Middle column: True state $\mathbf{X}$ in black and forward orbit of $\mathbf{x}_1$ under $\F$ in green.  Note green trajectory has been extended to $T_{max}$ five times that of training data.  Right column: True measurement noise $\mathbf{N}$ in blue, learned measurement noise $\mathbf{\hat{N}}$ in red, and error in estimation of measurement noise as dashed black line.}
\label{fig:navier_stokes_learned}
\end{figure}

We generate data by direct numerical simulation of the two-dimensional Navier-Stokes equations using the immersed boundary projection method~\cite{taira:07ibfs,taira:fastIBPM}, resulting in $151$ snapshots in time with spatial resolution of $199\times 449$.  
As in~\cite{Brunton2016}, we extract the time series of the first two proper orthogonal decomposition (POD) modes and the shift mode of Noack et al.~\cite{Noack2003jfm} as our clean training data; these modes are shown in Fig.~\ref{fig:NS_modes}.  
We add noise to the data following projection onto the low-dimensional subspace.  
The mean $L^2$ errors for the measurement noise approximation are $1.019$e$-4$ and $0.0504$ for the cases of 0\% and 1\% noise, respectively.  
We do not compute error metrics for vector field accuracy since the true vector field is unknown.  
However, the qualitative behavior of observations and model predictions match the training data, shown in Fig.~\ref{fig:navier_stokes_learned}. 

\subsection{Double pendulum}
In all of the example investigated so far, the true equations of motion have been simple polynomials in the state, and would therefore be easily represented with a sparse regression method such as the sparse identification of nonlinear dynamics (SINDy)~\cite{Brunton2016}.  
The utility of a neural network for approximating the vector field becomes more clear when we consider dynamics that are not easily represented by a standard library of elementary functions.  
The double pendulum is a classic mechanical system exhibiting chaos, with dynamics that would are  challenging for a library method, although previous work has recovered the Hamiltonian via genetic programming~\cite{Schmidt2009science}.  
The double pendulum may be modeled by the following equations of motion in terms of the two angles $\theta_1$ and $\theta_2$ of the respective pendula from the vertical axis and their conjugate momenta $p_1$ and $p_2$:
\begin{equation}
\begin{aligned}
\dot{\theta}_1 &= \dfrac{l_2p_1 - l_1p_2\cos (\theta_1 - \theta_2)}{l_1^2l_2 \left( m_1+m_2\sin^2 (\theta_1-\theta_2) \right)}\\
\dot{\theta}_2 &= \dfrac{-m_2l_2p_1 \cos(\theta_1 - \theta_2) + (m_1+m_2) l_1p_2}{m_2l_1l_2^2 \left( m_1+m_2\sin^2 (\theta_1-\theta_2) \right)}\\
\dot{p}_1 &= -(m_1+m_2)gl_1\sin(\theta_1) - C_1 + C_2\sin(2(\theta_1 - \theta_2)) \\
\dot{p}_2 &= -m_2gl_2 \sin(\theta_2) + C_1 - C_2 \sin(2(\theta_1 - \theta_2)),
\end{aligned} \label{eq:double_pendulum}
\end{equation}
where
\begin{equation}
\begin{aligned}
C_1 &= \dfrac{p_1p_2\sin(\theta_1 - \theta_2)}{l_1l_2 (m_1+m_2 \sin^2(\theta_1-\theta_2))}\\
C_2 &= \dfrac{m_2l_2^2p_1^2 + (m_1+m_2)l_1^2p_2^2 - 2m_2l_1l_2p_1p_2\cos (\theta_1 - \theta_2)}{2l_1^2l_2^2 (m_1+m_2 \sin^2(\theta_1-\theta_2))^2},
\end{aligned} \label{eq:double_pendulum_2}
\end{equation}
and $l_1=l_2=1$ are the lengths of the upper and lower arms of the pendulum, $m_1=m_2=1$ the respective point masses, and $g=10$ is the acceleration due to gravity.  
Numerical solutions to \eqref{eq:double_pendulum} are obtained using a symplectic integrator starting from the initial condition ${(\theta_1, \theta_2, p_1, p_2) = (1,0,0,0)}$ from $t=0$ to $t=50$ with a timestep of $\Delta t = 0.01$.  
A symplectic or variational integrator is required to ensure that energy is conserved along the trajectory~\cite{Yoshida1990pla,Marsden2001dmvi}.  
It is important to note that this initial condition represents a low-energy trajectory of the double pendulum with non-chaotic dynamics existing on a bounded region of phase space.  
Neither pendulum arm makes a full revolution over the vertical axis.  
For higher energies, the method presented in this paper did not yield satisfying results.

\begin{figure}
\centering
\includegraphics[width=.95\textwidth]{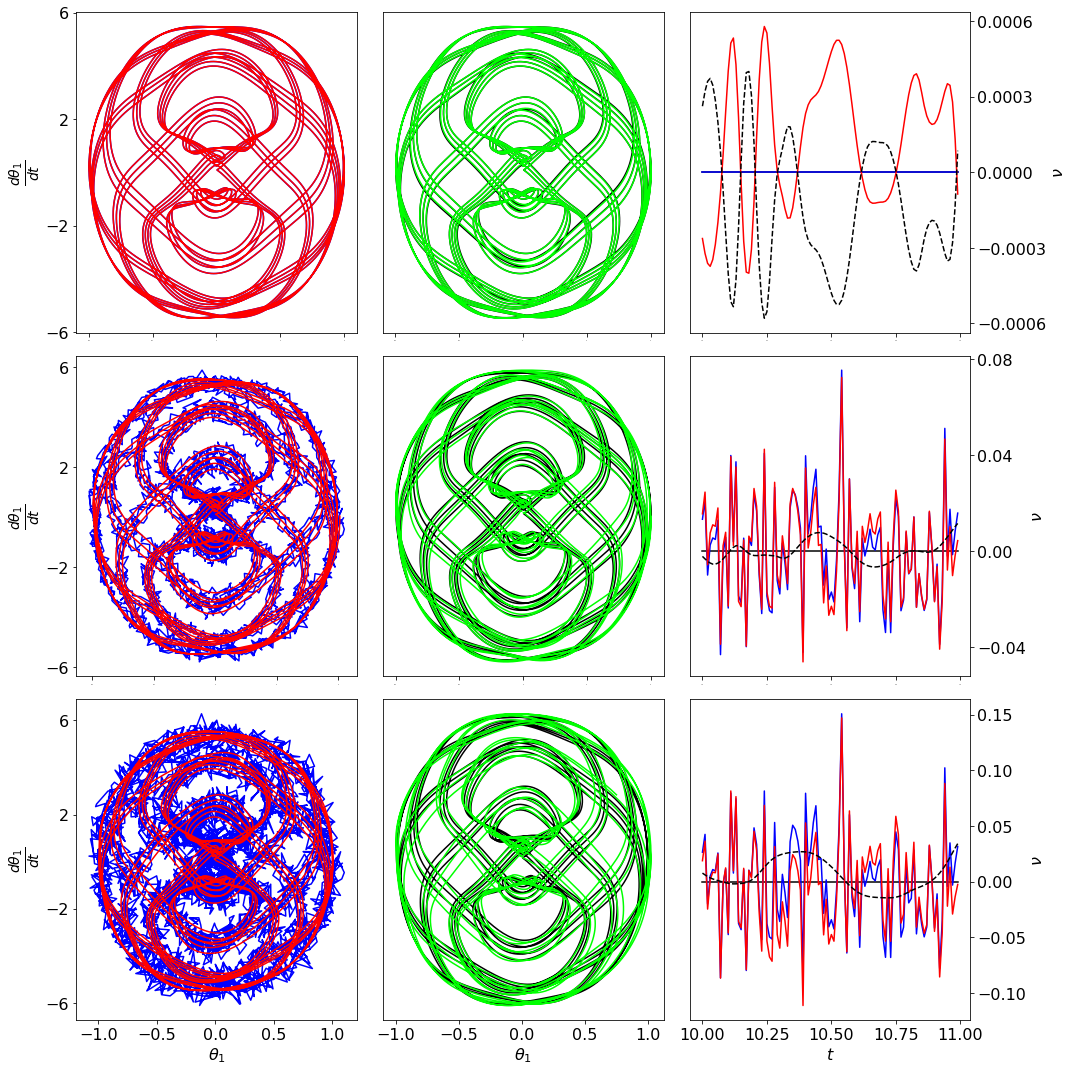}
\vspace{-.2in}
\caption{Results for double pendulum with increasing magnitudes of measurement noise.  Top row: Clean data.  Second row: 5\% Gaussian noise.  Bottom row: 10\% Gaussian noise.  Left column: Observations $\mathbf{Y}$ in blue and learned state $\mathbf{Y}-\mathbf{\hat{N}}$ in red.  Middle column: True state $\mathbf{X}$ in black and forward orbit of $\mathbf{x}_1$ under $\F$ in green.  Right column: True measurement noise $\mathbf{N}$ in blue, learned measurement noise $\mathbf{\hat{N}}$ in red, and error in estimation of measurement noise as dashed black line.}
\label{fig:double_pendulum_learned}
\vspace{-.1in}
\end{figure}

We construct a data-driven approximation to the vector field $f$ in \eqref{eq:double_pendulum} using a neural network for the vector field with five hidden layers, each containing $64$ nodes, embedded in a four-step Runge-Kutta scheme to approximate the discrete-time flow map $F$.  
Artificial measurement noise is added to the trajectory with magnitudes up to $10\%$.  Examples of training data, noise estimates, and numerical solutions to the learned dynamics are shown in Fig.~\ref{fig:double_pendulum_learned} for noise levels of $0$, $1$, and $10$ percent of the standard deviation of the dataset.  
Summary error measures for the learned dynamics and measurement noise of the double pendulum are shown in table \ref{tab:DP_error}.  
In all cases, it can be seen that the error is effectively separated from the training data, resulting in accurate model predictions, even for relatively large noise magnitudes.  

\begin{table}[b]
  \vspace{-.1in}
  \begin{center}
    \caption{Error for double pendulum with varying noise.}
    \label{tab:DP_error}
    \begin{tabular}{|l|c|c|c|c|}
      \hline
      \% Noise & $0$ & $1$ & $5$ & $10$ \\
      \hline
      $E_N$ & $5.823e-7$ & $5.835e-5$ & $1.200e-3$ & $3.399e-3$  \\
      $E_{f}$ & $5.951e-3$ & $6.192e-3$ & $8.575e-3$ & $1.444e-2$ \\
      \hline
    \end{tabular}
  \end{center}
  \vspace{-.1in}
\end{table}

\section{Cautionary remarks on neural networks and overfitting} \label{caution}
Neural networks are fundamentally interpolation methods~\cite{Mallat2016prsa} with high-dimensional parameter spaces that allow them to represent arbitrarily complex functions~\cite{hornik1991approximation}.  
The large number of free parameters required for arbitrary function fitting also creates the risk of overfitting, necessitating significant volumes of rich training data.  
Many successful applications of neural networks employ regularization techniques such as $L^2$ regularization, dropout \cite{srivastava2014dropout}, and early stopping \cite{caruana2001overfitting} to help prevent overfitting.  
In computer vision, data augmentation through preprocessing, such as random rotations and image flipping, helps prevent overfitting and allows single labeled examples of training data to be reused without redundancy.  
Recent innovations in the use of neural networks for dynamical systems forecasting have included regularization of the network Jacobian \cite{pan2018long}, but data augmentation does not have a clear analog in dynamical systems. 
This section will explore key limitations and highlight pitfalls of training neural networks for dynamical systems.  

Section \ref{results} demonstrated the ability of our proposed method to accurately represent dynamics from limited and noisy time-series data.  
In particular, the dynamics shown in the Lorenz equation, fluid flow, and double pendulum examples all evolved on an attractor, which was densely sampled in the training data.  
In each case, trajectories integrated along the learned dynamics remain on the attractor for long times.  
This indicates that our neural network is faithfully interpolating the vector field near the attractor.  
Fortunately, real-world data will often be sampled from an attractor, providing sufficiently rich data to train a model that is valid nearby. 
However, care must be taken when extrapolating to new initial conditions or when making claims about the vector field based on models learned near the attractor.  
In particular, transient data from off of the attractor may be essential to train models that are robust to perturbations or that are effective for control, whereby the attractor is likely modified by actuation~\cite{Brunton2015amr}.  

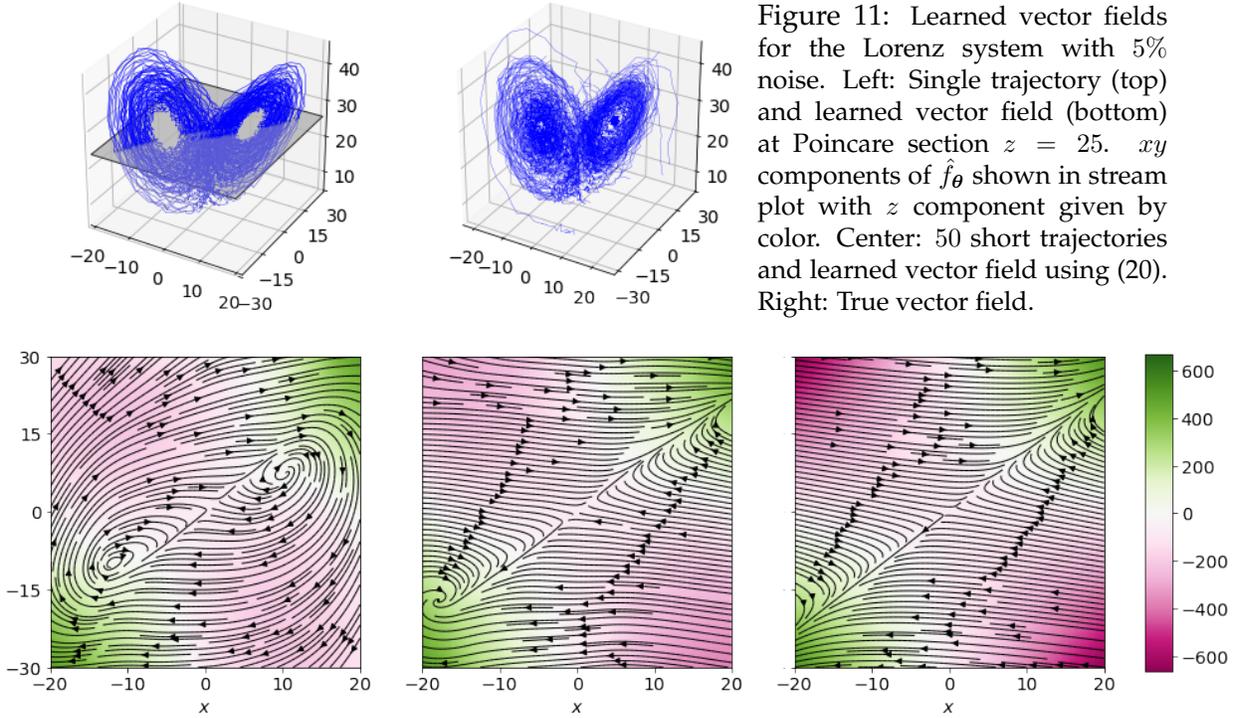
\begin{figure}
\centering
\includegraphics[width=\textwidth]{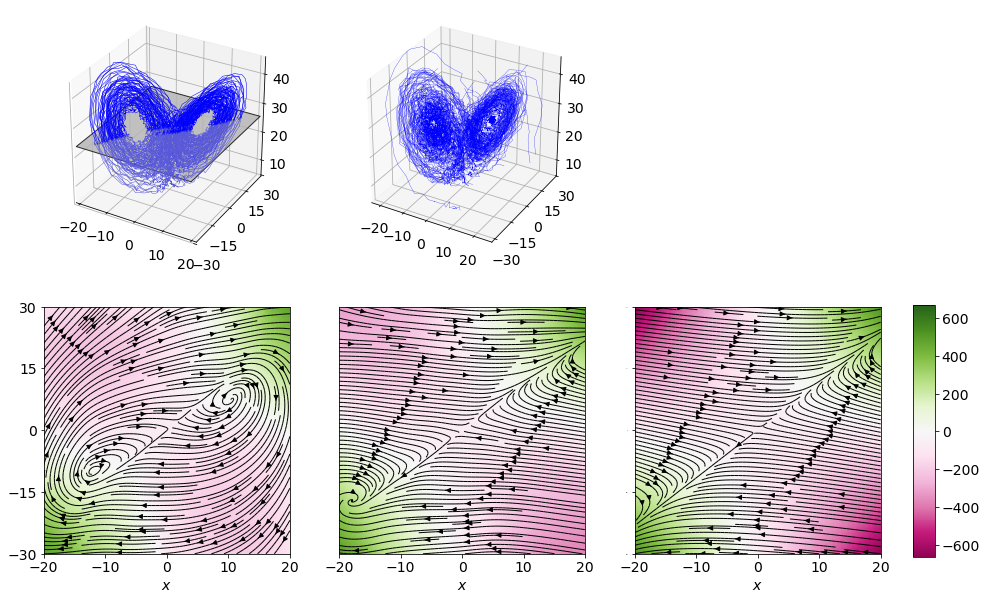}
\begin{picture}(0,0)
\put(55,230){\parbox{.33\linewidth}{
\noindent \caption{\small Learned vector fields for the Lorenz system with $5\%$ noise.  Left:  Single trajectory (top) and learned vector field (bottom) at Poincare section $z=25$.  $xy$ components of $\f$ shown in stream plot with $z$ component given by color.  Center:  $50$ short trajectories and learned vector field using \eqref{eq:multi_cost}.  Right:  True vector field.} \label{fig:lorenz_vector_field_slices}
}}
\end{picture}
\vspace{-.4in}
\end{figure}

By analyzing a known dynamical system such as the Lorenz equations, we are able to quantify the performance of the method in approximating the vector field when given only data on the attractor, or when given many trajectories containing transients.  
To do so, we extend our previous method to fit multiple datasets of observations to the same dynamical system.  
Given a set of $p$ trajectories, $\{\mathbf{Y}_k\}_{k=1}^p$, all having the same underlying dynamics, we adapt the cost function given by \eqref{eq:full_cost} to 
\begin{equation}
\mathcal{L}_{\text{multi}} \left(\btheta, \{\mathbf{\hat{N}}_k\}_{k=1}^p, \{\mathbf{Y}_k\}_{k=1}^p \right) = \sum_{k=1}^p \mathcal{L} (\btheta, \mathbf{\hat{N}}_k, \mathbf{Y}_k).
\label{eq:multi_cost}
\end{equation}
Based on \eqref{eq:multi_cost}, we compare the accuracy of the learned vector field from datasets of comparable size, obtained either from a single trajectory or from many short trajectories.  
Figure \ref{fig:lorenz_vector_field_slices} shows the data and learned vector field on the plane $z=25$ for the Lorenz system trained either from a single trajectory of length $m = 10000$, or from $50$ individual trajectories of length $200$, each from a random initial condition off the attractor.  The exact vector field is  shown for comparison.  
Unsurprisingly, the long trajectory does not result in accurate interpolation of the vector field off the attractor, while training from many trajectories with transients results in a more accurate model.

\begin{figure}
\vspace{-.3in}
\centering
\includegraphics[width=.95\textwidth]{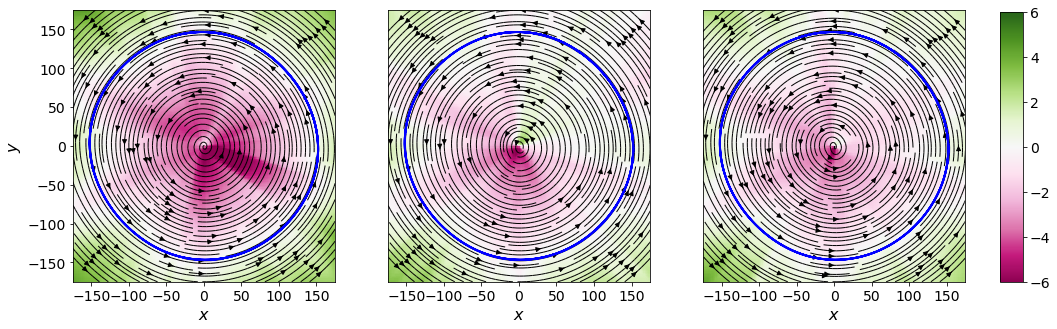}
\vspace{-.15in}
\caption{Learned vector fields for the reduced order system from flow around a cylinder with different sizes of neural network.  Left: 3 rows of 64 nodes each.  Center: 3 rows of 128 nodes.  Right:  5 rows of 64 nodes.  Mean field model exhibits radial symmetry.}
\label{fig:NS_vector_field}
\vspace{-.1in}
\end{figure}

Neural networks with different parameterizations will also result in varying behavior away from training data.  
To illustrate this, we consider three parameterizations of $f$ for the fluid flow dataset, using 3 hidden layers of size 64 or 128, as well as 5 hidden layers of size 64.  
The resulting vector fields along the $z=0$ plane are shown in Fig.~\ref{fig:NS_vector_field}.  
The limit cycle in the training data is shown in blue to indicate regions in the domain near training data.  
Each of the three parameterizations accurately models the forward orbit of the initial condition from the training data and converges to the correct limit cycle.  
However, all networks fail to identify a fundamental radial symmetry present in the problem, indicating that a single test trajectory is insufficient to enable forecasting for general time series.  
For data on the attractor, though, these specific parameterizations may be sufficient.

\section{Discussion} \label{discussion}
In this work we have presented a machine learning technique capable of representing the vector field of a continuous dynamical system from limited, noisy measurements, possibly with irregular time sampling.  
By explicitly considering noise in the optimization, we are able to construct accurate forecasting models from datasets corrupted by considerable levels of noise and to separate the measurement error from the underlying state.  
Our methodology constructs discrete time-steppers using a neural network to represent the underlying vector field, embedded in a classical Runge-Kutta scheme, enabling seamless learning from datasets unevenly spaced in time.  
The constrained learning architecture exploits structure in order to provide a robust model discovery and forecasting mathematical framework. 
Indeed, by constraining learning, the model discovery effort is focused on discovering the dynamics unencumbered by noise and corrupt measurements.

Using a neural network to interpolate the underlying vector field of an unknown dynamical system enables flexible learning of dynamics without any prior assumptions on the form of the vector field.  
This is in contrast to library-based methods, which require the vector field to lie in the span of a pre-determined set of basis functions, although using a neural network does forfeit interpretability.  
Both approaches have utility and may be thought of as complementary in terms of complexity and interpretability. 

The combination of neural networks and numerical time-stepping schemes suggests a number of high-priority research directions in system identification and data-driven forecasting.  
Future extensions of this work include considering systems with process noise, a more rigorous analysis of the specific method for interpolating $f$, including time delay coordinates to accommodate latent variables~\cite{Brunton2017natcomm}, and generalizing the method to identify partial differential equations. 
Rapid advances in hardware and the ease of writing software for deep learning will enable these innovations through fast turnover in developing and testing methods.  
To facilitate this future research, and in the interest of reproducible research, all code used for this paper has been made publicly available on GitHub at https://github.com/snagcliffs/RKNN.

\section*{Acknowledgments}
We acknowledge generous funding from the Army Research Office (ARO W911NF-17-1-0306 and W911NF-17-1-0422) and the Air Force Office of Scientific Research (AFOSR FA9550-18-1-0200).  

\small
\begin{spacing}{.5}
\bibliographystyle{plain}
\bibliography{nn_ode_bibliography}
\end{spacing}

\normalsize
\newpage
\section*{Appendix A: Expected error and structure of loss function}\label{expected_loss}
Our loss function in \eqref{eq:full_cost} evaluates the accuracy of a data-driven model against measured data using pairs of data separated by $i \in [1,q]$ timesteps.  For larger $i$, we expect errors to accumulate, resulting in less accurate predictions.  
Therefore, we use an exponential weight $\omega_i$ to discount larger $i$.  Here we estimate the error for an $i$ timestep prediction and show that under a relaxed set of assumptions, an upper bound is exponential in $i$, justifying our choice of $\omega_i$.  Let $\boldsymbol{\eta}_{j,i}$ be the error in approximating $\mathbf{y}_{j+i}$ from $\mathbf{y}_j$.  Then,
\begin{equation}
\boldsymbol{\eta}_{j,i} = \left(\F^i (\mathbf{y}_j - \boldsymbol{\hat{\nu}}_j ) + \boldsymbol{\hat{\nu}}_{j+1}\right) - \mathbf{y}_{j+i} .
\end{equation}
We are interesting in obtaining plausible estimates for the rate of growth of $\boldsymbol{\eta}_{j,i}$ as $i$ grows from $1$ to $q$.  Let $\boldsymbol{\epsilon}_j$ and $G$ denote the error in approximating the measurement noise and flow map respectively: 
\begin{equation}
\begin{aligned}
\boldsymbol{\hat{\nu}}_j - \boldsymbol{\nu}_j = \boldsymbol{\epsilon}_j\\
\F - F = G.
\end{aligned}
\end{equation}
We restrict our attention to a domain $\mathcal{D} \subset \R^n$ and let $\alpha = \sup_{\mathbf{x} \in \mathcal{D}} \|\mathbf{D}F(\mathbf{x})\|$ where $\mathbf{D}F$ is the Jacobian, $\tilde{\alpha} = \sup_{\mathbf{x} \in \mathcal{D}} \|\mathbf{D}\F(\mathbf{x})\|$,  and $\zeta = \sup_{\mathbf{x} \in \mathcal{D}} \|G(\mathbf{x})\|$.  We assume that the data-driven model is sufficiently accurate that $\zeta \ll 1$ and $\|\boldsymbol{\epsilon}_j\| < \mu \ll 1$.  The error in predicting $\mathbf{y}_{j+i}$ from $\mathbf{y}_j$ is,
\begin{equation}
\begin{aligned}
\boldsymbol{\eta}_{j,i} &= (\F^i (\mathbf{y}_j - \boldsymbol{\hat{\nu}}_j) + \boldsymbol{\hat{\nu}}_{j+1}) - \mathbf{y}_{j+1}\\
&= (\F^i (\mathbf{y}_j - \boldsymbol{\nu}_j -\boldsymbol{\epsilon}_j) + \boldsymbol{\nu}_{j+i} + \boldsymbol{\epsilon}_{j+i})  - \mathbf{y}_{j+i}\\
&= \F^i(\mathbf{x}_j -\boldsymbol{\epsilon}_j)  - \mathbf{x}_{j+i} + \boldsymbol{\epsilon}_{j+1}\\
&= \F^i(\mathbf{x}_j) - \mathbf{x}_{j+i} + \boldsymbol{\epsilon}_{j+1} + \mathcal{O}\left( \tilde{\alpha}^i \|\boldsymbol{\epsilon_j}\| \right)\\
&= (F+G)^i(\mathbf{x}_j) - F^i(\mathbf{x}_j) + \mathcal{O}\left( (\tilde{\alpha}^i +1)\mu \right).
\end{aligned}
\end{equation}
Focusing on the term $(F+G)^i(\mathbf{x}_j)$, let $\boldsymbol{\delta}_k = G((F+G)^{k-1}(\mathbf{x}_i))$ be the error in the data-driven flow map evaluated at $\F^{k-1}(\mathbf{x}_j)$. Then $\|\boldsymbol{\delta}_k\| \leq \zeta$ and,
\begin{equation}
\begin{aligned}
(F+G)(\mathbf{x}) &= F(\mathbf{x}) + \boldsymbol{\delta}_1 \\
(F+G)^2(\mathbf{x}) &= (F+G)(F(\mathbf{x}) + \boldsymbol{\delta}_1) \\
&\approx F^2(\mathbf{x}) + \mathbf{D}F(F(\mathbf{x}))\boldsymbol{\delta}_1 + \boldsymbol{\delta}_2 + \mathcal{O}(\zeta^2)\\
(F+G)^3(\mathbf{x}) &\approx (F+G)(F^2(\mathbf{x}) + \mathbf{D}F(F(\mathbf{x}))\boldsymbol{\delta}_1 + \boldsymbol{\delta}_2) \\
&= F^3(\mathbf{x}) + \mathbf{D}F(F^2(\mathbf{x}))\mathbf{D}F(F(\mathbf{x}))\boldsymbol{\delta}_1 + \mathbf{D}F(F^2(\mathbf{x}))\boldsymbol{\delta}_2 +\boldsymbol{\delta}_3+ \mathcal{O}(\zeta^2).
\end{aligned}
\end{equation}
Continuing to higher powers of $(F+G)$ we find,
\begin{equation}
\begin{aligned}
(F+G)^i(\mathbf{x})-F(\mathbf{x}) &= \sum_{k=1}^i \left( \prod_{l=1}^{i-k} \mathbf{D}F(F^{l}(\mathbf{x})) \right) \boldsymbol{\delta}_k \\
\|(F+G)^i(\mathbf{x})-F(\mathbf{x})\| &\leq \sum_{k=1}^i \left\|\left( \prod_{l=1}^{i-k} \mathbf{D}F(F^{l}(\mathbf{x})) \right) \boldsymbol{\delta}_k\right\| \\
&\leq \sum_{k=0}^{i-1} \zeta \alpha^{k} \\
&=\frac{\zeta(\alpha^i-1)}{\alpha-1},
\end{aligned}
\end{equation}
for $\alpha \neq 1$, and $(F+G)^i(\mathbf{x})-F(\mathbf{x}) \leq i\zeta$ otherwise.  Having ignored all quadratic terms in $\mu$ and $\zeta$, we find that an upper bound $\boldsymbol{\eta}_{j,i}$ is given by,
\begin{equation}
\| \boldsymbol{\eta}_{j,i} \| \leq \frac{\zeta(\alpha^i-1)}{\alpha-1} + \left( \tilde{\alpha}^i+1 \right) \mu , \label{eq:error_bound}
\end{equation}
where the norm used is the same as in the definitions of $\alpha$, $\tilde{\alpha}$, $\zeta$, and $\mu$. This expression is exponential in $i$. A similar expression may be derived for $i<0$ by following the same steps using the inverse flow maps.  Therefore we expect error to grow exponentially in the number of forward and backward timesteps.  In practice it would be very difficult to precisely estimate the quantity given in \eqref{eq:error_bound}, so we assume $\boldsymbol{\eta}_i = \boldsymbol{\eta}_0 \rho^i$ for some $\rho > 1$.  In particular, for chaotic systems and longer timesteps, one would use larger values of $\rho$, resulting in a more aggressive exponential discount in $\omega_i$. 

The accumulation of error through multiple iterations of the flow map informs how we weight our loss function.  In many machine learning tasks the canonical choice for loss function is mean square error, which is derived from the maximum likelihood estimate given an assumption of Gaussian distributed errors with identity covariance matrix.  For this work, we use the $L^2$ metric for our loss function and weight errors for $i$-step predictions with exponentially decreasing magnitude $\omega_i = \omega_0 \rho^{-i}$. 
This assumes Gaussian distributed error for all predictions, which is naive.  However, any standard metric would fall short of capturing the true error of our computational framework for any reasonably complex dynamical system.  Future work is required to investigate this more carefully.

\end{document}